\documentclass{l4dc2025}

\title[Neural Contraction Metrics with Formal Guarantees for Discrete-Time Systems]{
Neural Contraction Metrics with Formal Guarantees\\for Discrete-Time Nonlinear Dynamical Systems}
\usepackage{times}
\usepackage{algorithm}
\usepackage[noend]{algorithmic}
\usepackage[utf8]{inputenc}
\usepackage{booktabs}
\usepackage{multirow}
\usepackage{graphicx}
\usepackage{caption} 
\usepackage{todonotes}
\usepackage{subcaption} 
\usepackage{enumitem}
\usepackage{titlesec}
\titlespacing*{\section}{0pt}{1.7ex plus 0.5ex minus 0.2ex}{1.2ex plus 0.2ex}
\titlespacing*{\subsection}{0pt}{1.4ex plus 0.4ex minus 0.2ex}{1ex plus 0.2ex}
\titlespacing*{\subsubsection}{0pt}{1.2ex plus 0.3ex minus 0.2ex}{0.8ex plus 0.2ex}
\newcommand{\norm}[1]{\lVert{#1}\rVert}

\coltauthor{
\Name{Haoyu Li} \nametag{\footnotemark[1]\footnotetext[1]{Equal Contribution.}}\Email{haoyuli5@illinois.edu}\\
\Name{Xiangru Zhong}\nametag{\footnotemark[1]}\Email{xiangru4@illinois.edu}\\
\Name{Bin Hu} \Email{binhu7@illinois.edu}\\
\Name{Huan Zhang} \Email{huan@huan-zhang.com}\\
\addr University of Illinois Urbana-Champaign (UIUC), Department of Electrical and Computer Engineering (ECE) and Department of Computer Science (CS)
}

\begin{document}

\maketitle

\vspace{-1em}
\begin{abstract}%
Contraction metrics are crucial in control theory because they provide a powerful framework for analyzing stability, robustness, and convergence of various dynamical systems. However, identifying these metrics for complex nonlinear systems remains an open challenge due to the lack of scalable and effective tools. This paper explores the approach of learning verifiable contraction metrics parametrized as neural networks (NNs) for discrete-time nonlinear dynamical systems. While prior works on formal verification of contraction metrics for general nonlinear systems have focused on convex optimization methods (e.g. linear matrix inequalities, etc) under the assumption of continuously differentiable dynamics, the growing prevalence of NN-based controllers, often utilizing ReLU activations, introduces challenges due to the non-smooth nature of the resulting closed-loop dynamics. To bridge this gap, we establish a new sufficient condition for establishing formal neural contraction metrics for general discrete-time nonlinear systems assuming only the continuity of the dynamics. We show that from a computational perspective, our sufficient condition can be efficiently verified using the state-of-the-art neural network verifier $\alpha,\!\beta$-CROWN, which scales up non-convex neural network verification via novel integration of symbolic linear bound propagation and branch-and-bound.
Built upon our analysis tool, we further
develop a learning method for synthesizing neural contraction metrics from sampled data.  Finally, our approach is validated through the successful synthesis and verification of NN contraction metrics for various nonlinear examples.
\end{abstract}




\section{Introduction}
Contraction theory is a prominent framework for analyzing the stability and safety of control systems. It examines how trajectories of a system remain close and converge to each other~\citep{lohmiller1998contraction, tran2018convergence,giesl2022review,manchester2017control,tsukamoto2021contraction,bullo2022contraction,dawson2023safe}. In contrast to Lyapunov theory, which assesses stability with respect to a specific equilibrium point~\citep{lyapunov1992general}, contraction analysis focuses on the existence of a Riemannian-type metric~$d$ in which the system is contracting~\citep{simpson2014contraction}. This approach certifies that the distance between any pair of trajectories diminishes over time, ensuring that all solutions converge toward each other regardless of their initial conditions. Contraction theory has been applied to analyze and synthesize controllers for both continuous-time~\citep{tsukamoto2020neural, manchester2017control,sun2021learning,tsukamoto2021contraction} and discrete-time systems~\citep{wei2021control, wei2022discrete}. However, many of these approaches are often heuristic-based methods and lack formal verification.

In this paper, we focus on synthesizing and formally verifying a neural network (NN) parametrized contraction metric that is valid over the largest possible domain for a given discrete-time nonlinear dynamical system. Although previous studies such as \citet{giesl2019construction}, \citet{giesl2024contraction}, \citet{pokkakkillath2024construction}, and \citet{10714396} employ numerical quadrature or learning-based approaches to synthesize contraction metrics, they share common limitations.
First, they all implicitly assume the smoothness of the system dynamics. However, with the advancements of deep learning, NN-based controllers have been popularized~\citep{chang2019neural, dai2021lyapunov, wu2023neural, yang2024lyapunov,wang2024actor}, and the use of ReLU/LeakyReLU activations~\citep{dai2021lyapunov, wu2023neural, yang2024lyapunov} just naturally introduce non-smoothness into the closed-loop dynamics. Consequently, traditional synthesis and verification techniques that rely on the Jacobian matrix of the dynamics are no longer applicable. Second, many existing works depend on convex optimization techniques such as linear matrix inequalities for the verification and learning of contraction metrics, which unavoidably increases the conservativeness of the resulting metrics.


To overcome these challenges, we devise a novel verification condition for discrete-time contraction analysis that depends neither on matrix inequalities nor on the Jacobian of the dynamics. 
Besides handling nonsmooth systems, this novel formulation also significantly reduces the complexity of verification conditions. It allows the reformulation of contraction certification into a verification problem within a bounded set on a computation graph, making it well-suited for advanced neural network verification tools like $\alpha,\!\beta$-CROWN~\citep{zhang2018efficient,xu2020automatic,xu2021fast,wang2021beta}, which has demonstrated its exceptional scalability in the verification of neural network controllers~\citep{everett2023drip,yang2024lyapunov}. In contrast to traditional SMT and MIP solvers, $\alpha,\!\beta$-CROWN is far more efficient for this task, as it leverages the neural network architecture and benefits from GPU acceleration. By combining our novel verification formulation with the efficient neural network verifier, we achieved the first formally verified contraction metrics for a state feedback system with an NN controller. To summarize, our main contributions include:


\begin{itemize}[wide, labelwidth=!, labelindent=0pt, itemsep=-3pt, topsep=3pt]
    \item We develop a theoretical framework for formal contraction certificates applicable to potentially \textbf{non-smooth} discrete-time systems, proposing a \textbf{novel verification condition} for contraction metrics that eliminates the need for both matrix inequalities and the Jacobian of the dynamics;
    \item Unlike previous work with verified contraction metrics~\citep{10714396}, which uses expensive verifiers like SMT, our new formulation allows the combination of advanced neural network verification tool $\alpha,\!\beta$-CROWN, enabling us to use \textbf{effective learning-based NN contraction metrics} while ensuring efficient \textbf{formal verification of contraction};
    \item Our experiments on several autonomous and NN-controlled systems demonstrate that our learned NN contraction metric can produce a large verified region of contraction, sometimes close to the upper bound.
    For the first time in literature, we present verified neural contraction metrics for a \textbf{neural network-controlled state-feedback pendulum system}.
\end{itemize}


\section{Preliminaries}\label{sec:prelim}
\paragraph{Notations}Our paper primarily focuses on a discrete time-invariant closed-loop control system
\begin{align}\label{eq:system}
    x_{k+1} &= g(x_k, \pi(x_k)) = f(x_k)
\end{align}
where $x_k$ is the system state, $g$ denotes the open-loop dynamics, $\pi$ denotes the policy, and $f$ denotes the closed-loop nonlinear dynamics. We are interested in two cases: i) $f$ involves neural networks, or ii) $f$ is just a general nonlinear system. We denote the solution of the system starting with initial condition $x=x_0$ to be $\phi(k, x_0)$, i.e. it satisfies $\phi(k+1,x_0) = f(\phi(k, x_0))$.  We denote the set of all positive semi-definite/positive definite  matrices by $\mathbb{S}_{+}^{n\times n}$ and $\mathbb{S}_{++}^{n\times n}$ respectively. For two matrices $A$ and $B$, we write $A \prec B$ or $A \preceq B$ to indicate that $B-A$ is positive definite or positive semi-definite, respectively. A function $d : \mathcal{X} \times \mathcal{X} \to \mathbb{R}_{\geq 0}$ is said to be a metric on $\mathcal{X}$ if the following conditions holds: (1) $d(x,x) = 0, \,  \forall x \in \mathcal{X}$; (2) $d(x,y) = d(y,x), \forall x, y \in \mathcal{X}$; (3) $d(x,z) \leq d(x,y) + d(y,z), \forall x, y, z \in \mathcal{X}$. We denote the norm induced by a matrix $P \in \mathbb{S}_{++}^{n \times n}$ as $\Vert\cdot\Vert_P$ such that $\Vert x\Vert_P = \sqrt{x^{\top}Px}$, and we denote $B(x;\epsilon)$ as the infinity norm ball around $x$ with radius $\epsilon$, i.e. $B(0;\epsilon) = \{y: \Vert y-x\Vert_\infty \leq \epsilon\}$.

\paragraph{Contraction Metrics}
This paper focuses on contraction metrics, which describe the asymptotic property of the differences between two solutions of a system. The formal definition goes as follows. 

\begin{definition}
    [Contraction]
    The system~\eqref{eq:system} is said to be \textbf{contracting} on a forward invariant set $\mathcal{X} \subset \mathbb{R}^n$ if there exists a metric $d: \mathcal{X} \times \mathcal{X} \to \mathbb{R}_{\geq 0}$ and a scalar $\rho\in (0, 1)$ such that
    \begin{align}\label{def:contraction}
        d(f(x), f(y)) \leq \rho d(x,y), \forall x, y \in \mathcal{X}
    \end{align}
    We call $d$ the corresponding \textbf{contraction metric} and $\rho$ the \textbf{contraction rate}.
\end{definition}

\noindent
Clearly, if a system is contracting with some metric $d$, the distance between any two trajectories of the system will converge to $0$ with respect to the metric $d$.
In the simplest setting, one can adopt $d(x,y)=\norm{x-y}_2$, and the contraction analysis 
just reduces to standard Lipschitz constant analysis \citep{fazlyab2019efficient,huang2021training,araujo2023a,wang2024scalability}. However, general nonlinear systems typically require using contraction metrics in more general forms.
If the closed loop dynamics $f$ on the right hand side of system~\eqref{eq:system} is continuously differentiable, the existence of a metric $d$ in which the system is contracting can be reduced to existence of a matrix-valued function as covered in~\citet{lohmiller1998contraction}. The formal statement goes as follows.

\begin{theorem}\label{thm:matrix_version}
    The system~\eqref{eq:system} is contracting in $\mathcal{X}$ if there exists a nonsingular matrix-valued function $\Theta: \mathcal{X} \to \mathbb{R}^{n\times n}$ and positive constants $\mu, \eta, \rho$ such that for all $x \in \mathcal{X}$ we have
    \begin{equation}
        \eta I \leq \Theta(x)^\top\Theta(x) \preceq \rho I, \quad F(x)^\top F(x) - I \preceq -\mu I
    \end{equation}
    where $F(x)= \Theta(f(x))\frac{\partial f}{\partial x}(x) \Theta^{-1}(x)$.
\end{theorem}
For the proof of this theorem, we refer the readers to~\citet{tran2018convergence}. This criterion has several drawbacks. First of all, the verification of matrix inequalities can be hard for verifiers. The most common way of verifying matrix inequalities relies on Sylvester's criterion, as used by~\citet{10714396}. As can be seen from theorem~\ref{thm:matrix_version}, we need to verify multiple semi-definite conditions, and each of them requires us to verify a determinant condition for each of the $2^{n} - 1$ principal minors. Therefore, the verification problem cannot scale when we consider higher-dimension systems. Secondly, this condition requires the system to be continuously differentiable as it requires the Jacobian of the dynamics. However, currently, many state-of-the-art results have controllers synthesized by ReLU networks~\citep{yang2013neural,wu2023neural}, which makes this theorem inapplicable.



\paragraph{The $\alpha,\beta$-CROWN Verifier}
$\alpha,\beta$-CROWN~\citep{zhang2018efficient,xu2020automatic,xu2021fast,wang2021beta,zhang2022general} is a state-of-the-art neural network verification tool that aims to rigorously verify properties of neural networks. For a box $\mathcal{B}$ and a general function $F$ (which could involve NNs), it aims to prove $F(x) \geq 0$ for all $x \in \mathcal{B}$. The $\alpha,\!\beta$-CROWN verifier implements an efficient linear bound propagation method and branch-and-bound. The verification problem is solved if this linear lower bound is greater than zero on the domain $\mathcal{B}$. Otherwise, it splits the domain into subdomains and further computes bounds for each subdomain until the condition is verified.  More generally, if the problem contains multiple constraints such as $(F_1(x) \geq 0) \vee \cdots \vee (F_n(x) \geq 0)$, $\alpha,\!\beta$-CROWN will compute a linear lower bound for each of $F_1, \cdots, F_n$, and claim successful verification of the problem if at least one of the lower bounds is greater than $0$. Besides regular neural networks, the $\alpha,\!\beta$-CROWN toolbox supports many operations crucial in control applications such as trigonometric and polynomial functions.
Compared to traditional SMT and MIP solvers, $\alpha,\!\beta$-CROWN offers better scalability since its core linear bound propagation procedure exploits the structure of NNs and computation graphs. This procedure does not rely on an external mathematical programming solver and can be efficiently accelerated on GPUs.


\section{Contraction for Non-Smooth Dynamics}\label{sec:theory}
In this section, we aim to establish a condition that avoids the reliance on the Jacobian of the dynamics $f$ and can be efficiently handled by $\alpha,\!\beta$-CROWN. To achieve this, we propose conditions that guarantee contraction of the closed-loop dynamics under a Riemannian metric only assuming the continuity of $f$ and mild conditions on the domain of interest $\mathcal{X}$. This approach bridges contraction theory with our practical verification methodology and enables us to efficiently learn and verify neural contraction metrics without the computational limitations imposed by matrix inequalities.

\begin{theorem}\label{thm:main}
    Assume that $\mathcal{X}$ is an open, connected, and forward invariant subset in $\mathbb{R}^n$, and the dynamics $f$ is continuous. Given $\epsilon > 0$ a positive threshold. Suppose there exists uniformly continuous $M(x) : \mathcal{X} \to \mathbb{S}^{n\times n}_{++}$ with a uniform lower bound $\mu I \preceq M(x)$ that satisfies 
    \begin{align}\label{verification}
         \sqrt{(f(x) - f(y))^\top M(f(x))(f(x)-f(y))} \leq \rho \sqrt{(x-y)^\top M(x)(x-y)} 
    \end{align}
    for all $x \in \mathcal{X}$ and $y \in B(x; \epsilon) \cap \mathcal{X}$ and some $0< \rho < 1$. For every pair of nonnegative real numbers $a \leq b$, we define the space of admissible curves with starting time $a$ and end time $b$ as
    \begin{align}
        \Gamma_{a,b} := \{\gamma: [a,b] \to \mathcal{X}, \text{$\gamma$ is continuous and piecewise regular} \}
    \end{align}
    where piecewise regular means that on each piece we have $\gamma \in C^1$ and $||\gamma'(t)|| \neq 0$ for all $t$. The length of an admissible curve $\gamma: [a,b] \to \mathcal{X}$ under the metric $M$ is defined as
    \begin{align}
        L(\gamma) := \int_a^b \sqrt{\gamma'(t)^\top  M(\gamma(t))\gamma'(t)} \, dt = \int_a^b \Vert \gamma'(t)\Vert_{M(\gamma(t))} \, dt
    \end{align}
    For each $x,y$, we then define $\Gamma_{x,y} := \bigcup_{a,b} \{\gamma \in \Gamma_{a,b}, \gamma(a) = x, \gamma(b) = y\}$.
     It is then guaranteed that the system~\eqref{eq:system} is contracting on $\mathcal{X}$ in the metric $d(x,y) = \inf_{\gamma \in \Gamma_{x,y}} L(\gamma)$. 
\end{theorem} 

\begin{lemma}
    For any continuous curve $\gamma: [a,b] \to \mathcal{X}$ that connects $x, y \in \mathcal{X}$, there exists an admissble curve $\gamma_s: [a,b] \to \mathcal{X}$ connecting $x, y$ and $||\gamma(t) - \gamma_s(t)|| \leq \epsilon$ for any $\epsilon > 0$ and any $t$.
\end{lemma}

\noindent
\textbf{Proof:} 
Fix $\epsilon > 0$. Since $\gamma$ is continuous, the image $\gamma([a,b])$ is compact. Since $\mathcal{X}$ is assumed to be open, for each $t \in [a,b]$, there exists an $\epsilon_t > 0$ such that $B(\gamma(t);\epsilon_t) \subset \mathcal{X}$. Now we consider the open cover of $\gamma([a,b])$ by $\{B(\gamma(t); \frac{\epsilon_t}{2})\}_t$. By compactness, it must admit a finite subcover $B(\gamma(t_k);\frac{\epsilon_{t_k}}{2})$. Now pick $r > 0$ such that $r < \
\min(\min_k{\frac{\epsilon_{t_k}}{2}}, \epsilon)$. As $\gamma$ is uniformly continuous on $[a,b]$, there exists a $\delta$ such that $|t -s| < \delta$ implies that $|\gamma(t) - \gamma(s)| < \frac{r}{4}$. Now we take a partition of $[a,b]$ into $N$ subintervals $\{[t_{i}, t_{i+1}]\}_{i=1}^N$ with $\gamma(t_1) = a$ and $\gamma(t_{N+1}) = b$ such that $\frac{\delta}{2} < |t_{i+1} -t_i| < \delta$ for any $i$. We now consider the curve that is a piecewise straight line on each subinterval $[t_i, t_{i+1}]$. More formally, with the partition at hand, we first construct $N+1$ points inductively as follows: We first choose $\phi_0 = \gamma(0)$. For $i \leq N$, we choose $\phi_i = \gamma(t_i)$ if $\gamma(t_i) \neq \phi_{i-1}$, and otherwise we chose $\phi_i$ in $B(\gamma(t_i);\frac{r}{16})$ such that $\phi_{i} \neq \phi_{i-1}$. For $i = N$ we additionally need to make sure that $\phi_{i} \neq \gamma(t_{N+1})$. And finally we pick $\phi_{N+1} = \gamma(t_{N+1})$. Now for $t \in [t_{i}, t_{i+1}]$ we define
\begin{align}
    \gamma_s(t) &= \phi_i + \frac{t-t_i}{t_{i+1} - t_i}(\phi_{i+1} - \phi_i)
\end{align}
Notice that the above construction makes sure that $\phi_{i+1} \neq \phi_i$ for all $i$. Therefore, on each piece we have $\gamma_s'(t) = \frac{\phi_{i+1} - \phi_i}{t_{i+1} - t_i} \neq 0$. We can now estimate that on each subinterval, we have
\begin{equation}
\begin{aligned}
    \Vert\gamma_s(t) - \gamma(t)\Vert &\leq \Vert\gamma_s(t) - \gamma(t_i)\Vert + \Vert\gamma(t) - \gamma(t_i)\Vert \\
                                &\leq \Vert \phi_{i+1} - \gamma(t_{i+1})\Vert + \Vert\gamma(t_{i+1}) - \gamma(t_i)\Vert + \Vert\gamma(t_i) - \phi_i\Vert + \Vert \phi_i - \gamma(t_i)\Vert + \frac{r}{4}\\
                                &<r \leq \epsilon
\end{aligned}
\end{equation}
Finally we check that $\gamma_s(t) \in \mathcal{X}$ for any $t$. As each $\gamma(t)$ lies in some $B(\gamma(t_k);\frac{\epsilon_{t_k}}{2})$, we have
\begin{align}
    \Vert\gamma_s(t) - \gamma(t_k)\Vert \leq \Vert\gamma_s(t) - \gamma(t)\Vert + \Vert\gamma(t) - \gamma(t_k)\Vert \leq r + \frac{\epsilon_{t_k}}{2} \leq \epsilon_{t_k}
\end{align}
So clearly by construction we have $\gamma_s(t) \in B(\gamma(t_k);\epsilon_{t_k}) \subset \mathcal{X}$. \hfill $\blacksquare$

\medskip
\noindent
\textbf{Proof of Thm~\ref{thm:main}:} For the proof that the $d$ is indeed a metric, we refer the readers to the text~\cite{lee2018introduction}. We shall now show that the verification condition~\eqref{verification} implies the contraction in the defined metric. For any $\gamma \in \Gamma_{x,y}$, since $\gamma$ is assumed to be regular on each piece, there exists a $\gamma_u: [c,d] \to \mathcal{X}$ such that $\gamma_u \in \Gamma_{x,y}$, $L(\gamma) = L(\gamma')$, and $\Vert\gamma_u'(t)\Vert = 1$ for almost every $t$ (this is arclength parametrization, see~\cite{lee2018introduction}). Now we consider the curve $f(\gamma_u)$. As $\mathcal{X}$ is forward invariant, $f(\gamma_u)$ also belongs to $\mathcal{X}$. Now the difficulty arises since $f(\gamma_u)$ might not be admissible, since we only assume $f$ is continuous. However, a small perturbation on $f(\gamma_u)$ is always possible to make it admissible without affecting its length too much. More formally, we fix a $\rho' \in (\rho,1)$. By the verification condition~\eqref{verification}, we know that for all  small enough $\delta_t$ we always have 
\begin{equation}\label{condition1}
\begin{aligned}
    (f(\gamma_u(t+\delta_t)) &- f(\gamma_u(t)))^\top M(f(\gamma_u(t+\delta_t)))(f(\gamma_u(t+\delta_t)) - f(\gamma_u(t))) \\ &\leq \rho^2 (\gamma_u(t+\delta_t)-\gamma_u(t))^\top M(\gamma_u(t+\delta_t))(\gamma_u(t+\delta_t)-\gamma_u(t))\\
    &< \rho'^2 (\gamma_u(t+\delta_t)-\gamma_u(t))^\top M(\gamma_u(t+\delta_t))(\gamma_u(t+\delta_t)-\gamma(t))
\end{aligned}
\end{equation}
Notice that per our assumption, for every fixed $t$ we can bound 
\begin{equation}
\begin{aligned}
    (\rho'^2-\rho^2) (\gamma_u(t+\delta_t)-\gamma_u(t))^\top M(\gamma_u(t+\delta_t))(\gamma_u(t+\delta_t)-\gamma_u(t)) &\geq \mu (\rho'^2 - \rho^2)\frac{\delta_t^2}{2}
\end{aligned}
\end{equation}
for all small $\delta_t$. Therefore, we see that after dividing by $\delta_t^2$, the margin of the strict inequality~\eqref{condition1} is at least $\frac{\mu}{2}(\rho'^2 -\rho^2)$. Now fix an $\epsilon > 0$, we know that $f(\gamma_u)$ is a continuous curve, so that we can construct an admissible curve $\gamma_p: [c,d] \to \mathcal{X}$ such that $\Vert f(\gamma_u(t)) - \gamma_p(t)\Vert \leq \epsilon$ for any $t$. Now we see that we have
\begin{equation}
\begin{aligned}
     |(f(\gamma_u(t+\delta_t)) - f(\gamma_u(t))) &- (\gamma_p(t+\delta_t) - \gamma_p(t))| \\ &\leq |f(\gamma_u(t+\delta_t)) - \gamma_p(t+\delta_t)| + |f(\gamma_u(t)) - \gamma_p(t)| \leq 2\epsilon.
\end{aligned}
\end{equation}
Now by~\eqref{condition1}, the inequality is strict when the contraction factor is $\rho'$. This gives us room to wiggle the curve $f(\gamma_u)$ to be $\gamma_p$, which clearly belongs to the set $\Gamma_{f(x),f(y)}$. When $\epsilon$ is chosen to be small enough, we shall still have the inequality condition~\ref{condition1} for any $t$ and small enough $\delta_t$, i.e.
\begin{equation}
\begin{aligned}
    \Vert\gamma_p(t+\delta_t) - \gamma_p(t)\Vert_{M(\gamma_p(t+\delta_t))}^2 \leq \rho'^2 \Vert\gamma_u(t+\delta_t) -\gamma_u(t)\Vert_{M(\gamma_u(t+\delta_t))}^2,
\end{aligned}
\end{equation}
as $M$ is uniformly continuous. Now send $\delta_t^2 \to 0$ we get 
\begin{align}
\Vert \gamma_p'(t)\Vert_{M(\gamma_p(t))} \leq \rho' \Vert \gamma_u'(t)\Vert_{M(\gamma_u(t))}.
\end{align}
This then implies 
\begin{align}
    d(f(x),f(y)) &\leq \int_c^d \Vert \gamma'_p(t)\Vert_{M(\gamma_p(t))} \, dt 
    \leq \rho'\int_c^d \Vert\gamma_u'(t)\Vert_{M(\gamma_u(t))} \, dt = \rho' L(\gamma)
\end{align}
Now, taking infimum over $\gamma$ on the right-hand side gives the desired contraction result. \hfill $\blacksquare$

\section{Methodology}\label{sec:method}
In this section, we present techniques to bridge our novel contraction theory in Section~\ref{sec:theory} and $\alpha,\beta$-CROWN to efficiently learn and verify neural contraction metrics.

\subsection{Invariant Set Finding}\label{subsection:Invariant Set Finding}
We first need to identify a forward invariant set 
$\mathcal{X}$ where the computation of contraction metrics will be performed. We adopt a method that demonstrates state-of-the-art results in certifying the region of exponential attraction to an equilibrium point in discrete-time systems~\citep{yang2024lyapunov}. We use the following theorem to find a forward invariant set $\mathcal{X}$ inside a box of interest $\mathcal{B}$:

\begin{theorem}
    Let $F(x) := V(f(x)) - (1-\kappa)V(x)$ where $\kappa > 0$. If the condition
    \begin{align}\label{veri:lyapunov_condition}
        (-F(x) \geq 0 \wedge f(x) \in \mathcal{B}) \vee (V(x) \geq \rho)
    \end{align}
    holds for any $x \in \mathcal{B}$, the closed loop dynamics~\eqref{eq:system} is Lyapunov stable with $V$ being a valid Lyapunov function and the set $\mathcal{X} := \{x: V(x) < \rho\} \cap \mathcal{B}$ is certified forward invariant set as well as an exponential region of attraction.
\end{theorem}

\noindent
The proof can be found in~\cite{yang2024lyapunov}. Eq. \eqref{veri:lyapunov_condition} is then converted into a loss function and trained with a counterexample-guided (CEGIS) approach to obtain the region of attraction.

\subsection{Formal Verification of Contraction Metrics}\label{formal_verify}
Given a candidate metric $M(x)$ (possibly a NN), following our theory in Sec.~\ref{sec:method}, we now derive the verification formulation such that a verifier like $\alpha,\!\beta$-CROWN can formally check if $M(x)$ is a contraction metric on the forward invariant domain $\mathcal{X}$. We introduce the auxiliary set $\mathcal{B}$ as $\alpha,\!\beta$-CROWN cannot handle implicitly defined domains like $\{V(x) < \rho\}$ directly. The verification condition for contraction can be formalized as follows.

\begin{theorem}
    Suppose we have already formally verified the forward invariance of the set $\mathcal{X} = \{x: V(x) < \rho\} \cap \mathcal{B}$. We define 
    \begin{align}\label{eq:G}
    G(x,\delta) := (f(x)-f(x+\delta))^\top M(f(x))(f(x)-f(x+\delta)) - \rho^2\delta^\top M(x)\delta.
    \end{align}
    If the condition
    \begin{align}\label{v:condition}
        (-G(x,\delta) \geq 0) \vee (x+\delta \not \in \mathcal{B}) \vee (V(x) \geq \rho) \vee (V(x+\delta) \geq \rho)
    \end{align}
    holds for any $(x,\delta) \in \mathcal{B} \times B(0;\epsilon)$ where $\epsilon > 0$, the system~\eqref{eq:system} is certifiably contracting on the domain $\mathcal{S}$ with metric $d$ and contraction rate $\rho' \in (\rho, 1)$.
\end{theorem}

\noindent
\textbf{Proof:} Condition~\eqref{v:condition} enforces equation~\eqref{eq:G} to hold for all $x \in \mathcal{X}$ and $y \in B(0;\epsilon)\cap \mathcal{X}$. The contraction is then implied by Theorem~\ref{thm:main}. \hfill $\blacksquare$

\smallskip
\noindent
This verification condition~\eqref{verification} allows us to verify contraction without relying on any Jacobian type condition or matrix inequality conditions, which enables the efficient verification by $\alpha,\!\beta$-CROWN. Each term of the condition~\eqref{verification} can be represented as constraints on general functions. We then leverage $\alpha,\!\beta$-CROWN's efficient bound propagation and branch-and-bound method to verify these constraints over the domain $\mathcal{B}$. While existing solvers like SMT solvers struggle with scalability, our approach is inherently more scalable and benefits from GPU acceleration as a byproduct, further enhancing its efficiency.

\subsection{Learning Neural Network based Contraction Metrics}\label{learning}

We now present a learning-based approach to synthesize the contraction metric $M(x)$. Given a forward invariant domain $\mathcal{X}$, we desired to synthesize a uniformly continuous $M: \mathcal{X} \to \mathbb{R}^{n\times n}$ that has a uniform positive lower bound. To ease the learning and verification, we enforce the definiteness condition and the lower bound condition by parametrizing $M$ to be
\begin{align}\label{M:param}
    M(x) := \mu I + R(x)^\top R(x)
\end{align}
where $\mu > 0$, and $R$ is a learnable function parametrized by a neural network. This ensures that $M$ is symmetric and $M(x) \geq \mu I$ for every $x \in \mathcal{X}$. The learning objective of $M(x)$ should be the satisfaction of the condition~\eqref{v:condition}, as demonstrated in the numerical condition below:
\begin{theorem}
    Except the condition $x+\delta \in \mathcal{B}$, the rest of the verification condition~\eqref{v:condition} is equivalent to the following numerical condition:
    \begin{align}\label{v:num_condition}
        L_{\text{violate}}(x,\delta; V, \rho) := \min(G(x,\delta), \rho - V(x), \rho - V(x+\delta)) \leq 0, \quad \forall x \in \mathcal{X}
    \end{align}
\end{theorem}

\smallskip
\noindent
As the condition $x + \delta \in \mathcal{B}$ can be easily verified posthoc, this numerical condition provides a simple training loss for the training of contraction metrics $M$, which is encoded in the function $G(x,\delta)$. We then adopt a counterexample-guided approach, where counterexamples of condition~\eqref{v:num_condition} are found and act as a dataset for the training. Instead of involving a sound verifier to generate the counterexamples as in~\cite{chang2019neural}, we utilized cheap projected gradient descent (PGD) attacks to find counterexamples rather than involving an inefficient verifier at this stage. We found that the cheap attacker will effectively find counterexamples without damaging the learning performance. Also, since it is unknown how large the region where the system can be certified to be contracting is, we gradually increase the size of the forward invariant set $\mathcal{X}$ during training. The training algorithm is shown in algorithm~\ref{alg:learn_contraction}. 


\begin{algorithm}
\small
\caption{Learning Contraction Metrics}
\label{alg:learn_contraction}
\begin{algorithmic}[1]
\REQUIRE metrics parameter $\theta$, dynamics $f$, region of interest $\mathcal{B}$, Lyapunov function $V$, levelset threshold $\rho$, learning rate $\eta$, PGD steps and stepsize $n$ and $\beta$, epochs $m$ 
\ENSURE A learned contraction metric, parameterized by $\theta$
\FOR{$\textit{scale }$ = 0.1, 0.2, $\cdots$, 1.0}
\STATE $\rho' = \textit{scale}\times \rho$ \hfill \textit{// Scale the threshold to adjust the contraction region}
\STATE Training dateset $\mathcal{D} = \emptyset$ \hfill \textit{// Initialize new empty dataset}
\FOR{$\textit{iter }$= 1, 2, $\cdots$}
\STATE $(\xi^i_{x}, \xi^i_{\delta}) \gets$ a batch of random PGD initializations
\FOR{$\textit{descent }$ = 1, 2, $\cdots$, $n$}
\STATE $L \gets L_{\text{violate}}(\xi^i_{x}, \xi^i_{\delta};V,\rho')$
\STATE $\xi^i_{x} \gets \text{Project}_{\mathcal{B}}(\xi^i_{x} + \beta\cdot\frac{\partial L}{\partial\xi^i_{x}})$ \hfill \textit{// PGD update for $x$} 
\STATE $\xi^i_{\delta} \gets \text{Project}_{B(0;\epsilon)}(\xi^i_{\delta} + \beta\cdot\frac{\partial L}{\partial \xi^i_{\delta}})$ \hfill \textit{// PGD update for $\delta$}
\ENDFOR
\STATE $\mathcal{D} \gets \mathcal{D} \cup \{(\xi^i_{x}, \xi^i_{\delta})\}$ \hfill \textit{// Update Dataset with counterexamples}
\FOR{$\textit{epoch } = 1, 2, \cdots, m$}
\STATE $L_{\text{train}} \gets \sum_{\mathcal{D}}  L_{\text{violate}}(\xi^j_x, \xi^j_\delta; V, \rho)$ \hfill \textit{// Summing over all the violations in the dataset}
\STATE $\theta \gets \theta - \eta \nabla_{\theta} L_{\text{train}}$ \hfill \textit{// Gradient Descent to optimize metric parameters}
\ENDFOR
\ENDFOR
\ENDFOR
\end{algorithmic}
\end{algorithm}

\noindent
\textbf{Remark:} We note that it is better not to include a strong regularization term on the parameters of the NN metric. Intuitively, suppose in addition that the forward invariant set $\mathcal{X}$ is convex and bounded, the contraction metric admits an upper bound $M(x) \preceq \eta I$, in addition to the lower bound $\mu I \preceq M(x)$. One can then prove that $\sqrt{\mu} \Vert x-y \Vert \leq d(x,y) \leq \sqrt{\eta} \Vert x-y \Vert$ and $\Vert \phi(k,x_0) - \phi(k,x_0')\Vert \leq \rho \sqrt{\frac{\eta}{\mu}}\Vert \phi(k-1,x_0) - \phi(k-1,x_0')\Vert$. The regularization encourages $|\eta - \mu|$ to be smaller, potentially enforcing a contraction condition with respect to the Euclidean norm. This condition is much stricter than contraction in metric $d$, and thus could be harmful to the learning.

\section{Experiments}\label{sec:exp}
To demonstrate the effectiveness of our method, we evaluate four different discrete-time systems that contain both smooth autonomous systems and a nonsmooth neural network state feedback system. All four systems are discretized with the explicit Euler method with a step size of 0.05. For each system, we first identify a Lyapunov function and a verified region of attraction using the method described in section~\ref{subsection:Invariant Set Finding}. Next, we train a contraction metric $M(x)$ following the procedure described in section~\ref{learning}. The local contraction condition~\eqref{v:condition} is verified using $\alpha,\!\beta$-CROWN with the specification described in~\ref{formal_verify}. For comparison, we also compute a constant contraction metric for each system by directly using $M(x) =  P$ provided by a quadratic Lyapunov function $V(x) = x^\top Px$ of the system. We shall demonstrate that by combining our training and verification, a large certified region of contraction can be achieved in all cases. In table~\ref{tab:overall}, we summarize the verified region of the learned NN contraction metric and its runtime.

\begin{table}[tb]
\centering
\small
\begin{tabular}{c|c|cc|cc|c}
\toprule
System            & ROA & \multicolumn{2}{c|}{NN metric (ours)}  & \multicolumn{2}{c|}{Constant metric (baseline)}  & Verification Runtime\\
                  & $\rho$ & $\rho$ ($\uparrow$)   & $r$ ($\uparrow$)   & $\rho$  ($\uparrow$)     & $r$ ($\uparrow$)     & NN Metric \\ \midrule
Van der Pol       & 0.119  & 0.119     & 100\%              & 0.04         & 33.7\%   &   109s         \\
polynomial        & 10.8   & 10.0      & 92.8\%             & 1.5          & 13.8\%   & 290s            \\
two-machine power & 0.086  & 0.086     & 100\%              & 0.027        & 31.3\%   & 62s \\
\midrule
inverted pendulum & 672    & 580       & 88.6\%             & invalid            & invalid         &5.8hrs         \\ \bottomrule
\end{tabular}
\caption{Verified regions of contraction of neural contraction metrics and verification runtime for different systems. These results demonstrate the effectiveness of the neural network contraction metrics, as well as the efficiency in terms of verification. The $\rho$ of ROA is the largest region we consider, and thus serves as an \textbf{upper bound} of what we can achieve. The ratio of the volume of each verified region to the volume of ROA is also presented, denoted as $r$. We can see that the NN metric achieves a significantly larger verified region compared to the baseline (illustrated in Fig.~\ref{fig:all_columns}). ``Invalid" means that the synthesized metric does not provide a nontrivial region of contraction.} 

\label{tab:overall}
\end{table}

\subsection{Nonlinear Autonomous Dynamical System}
\begin{figure}[t]

    \centering
    \textbf{NN Contraction Metric} \\[2pt]
    \begin{minipage}[t]{0.22\textwidth}
        \centering
        \includegraphics[width=\textwidth]{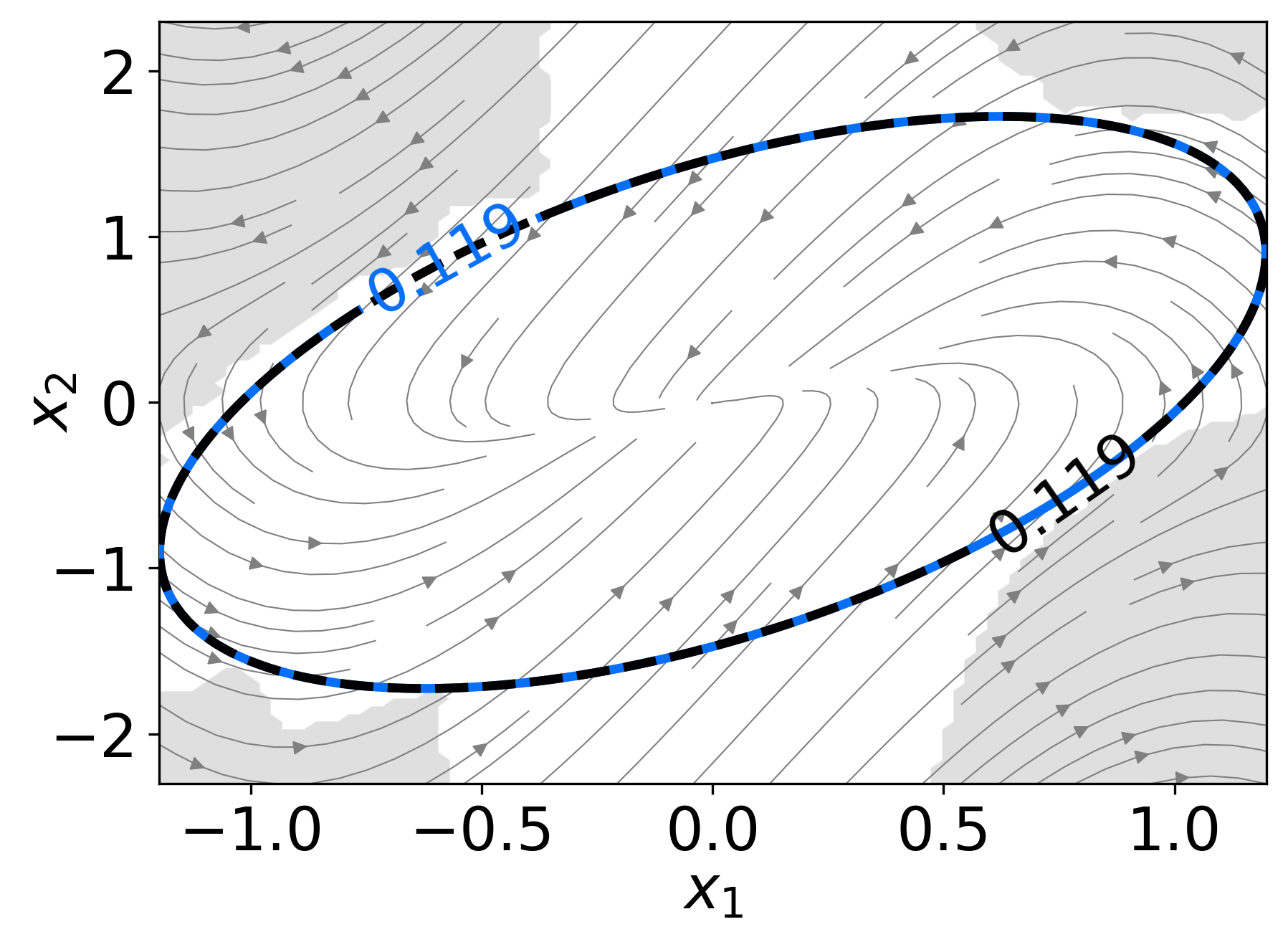}
    \end{minipage}
    \hfill
    \begin{minipage}[t]{0.22\textwidth}
        \centering
        \includegraphics[width=\textwidth]{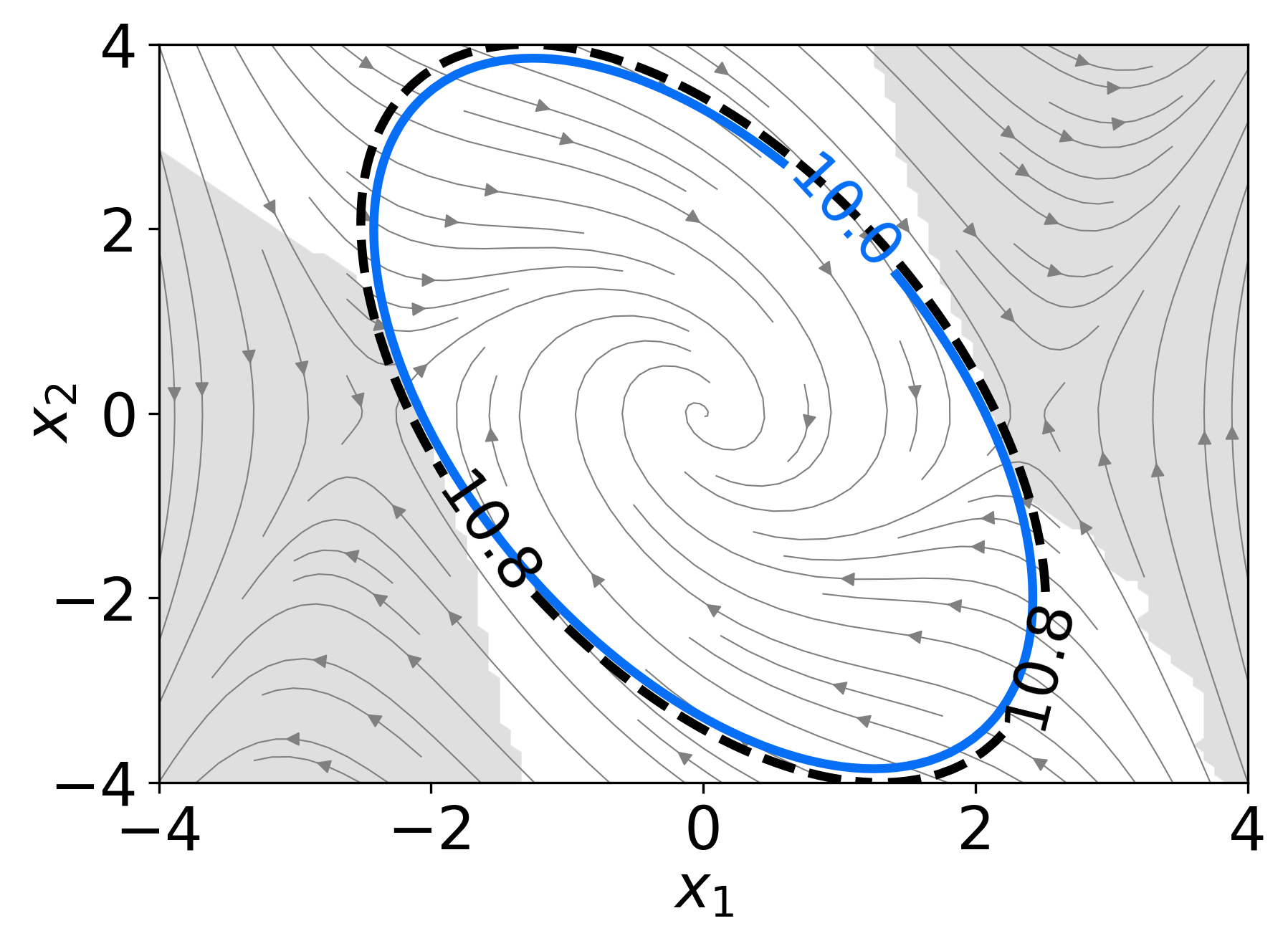}
    \end{minipage}
    \hfill
    \begin{minipage}[t]{0.22\textwidth}
        \centering
        \includegraphics[width=\textwidth]{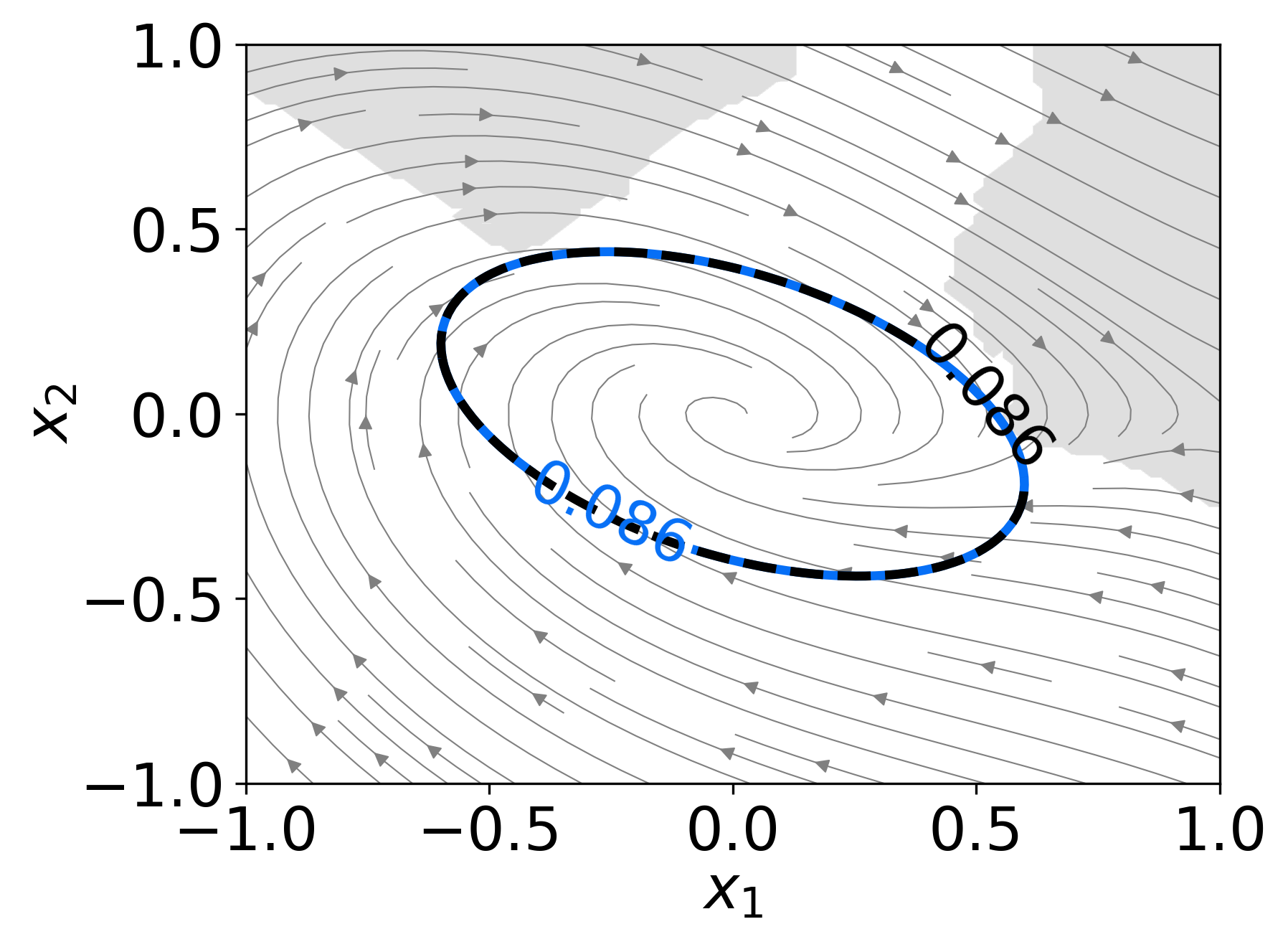}
    \end{minipage}
    \hfill
    \begin{minipage}[t]{0.22\textwidth}
        \centering
        \includegraphics[width=\textwidth]{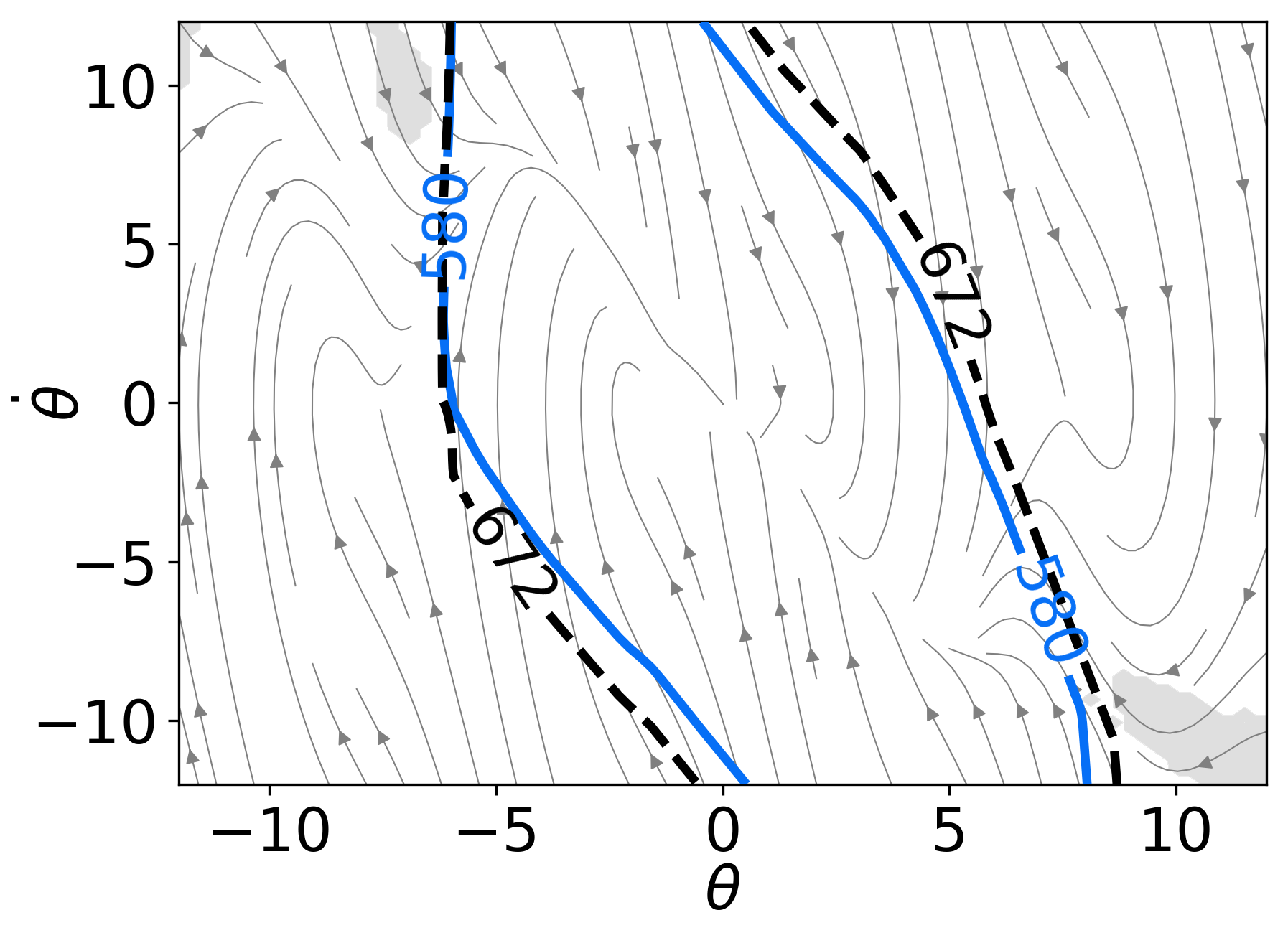}
    \end{minipage}


    \textbf{Constant Metric} \\[2pt]
    \begin{minipage}[t]{0.22\textwidth}
        \centering
        \includegraphics[width=\textwidth]{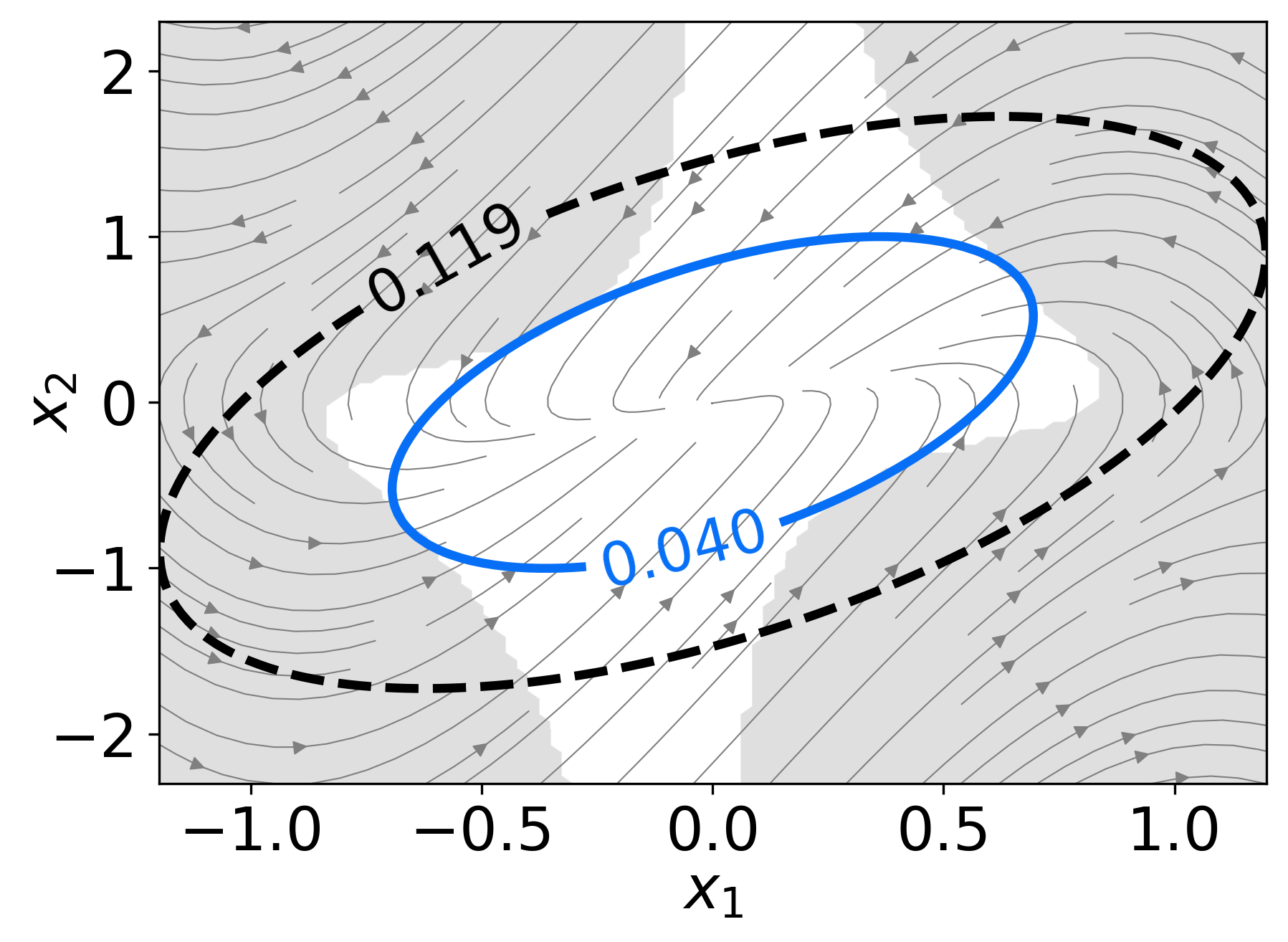}
        (a) Van der Pol
    \end{minipage}
    \hfill
    \begin{minipage}[t]{0.22\textwidth}
        \centering
        \includegraphics[width=\textwidth]{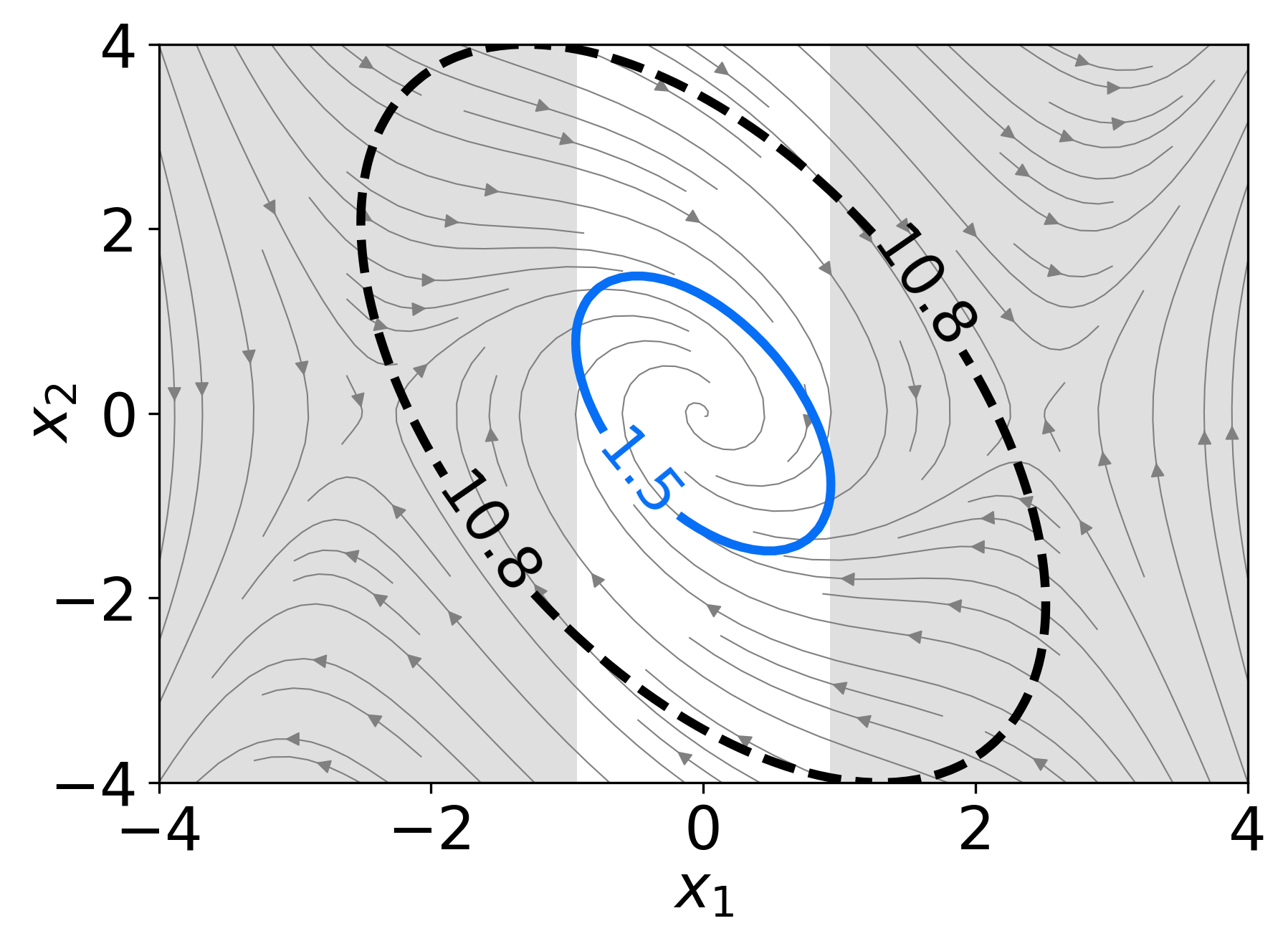}
        (b) Polynomial
    \end{minipage}
    \hfill
    \begin{minipage}[t]{0.22\textwidth}
        \centering
        \includegraphics[width=\textwidth]{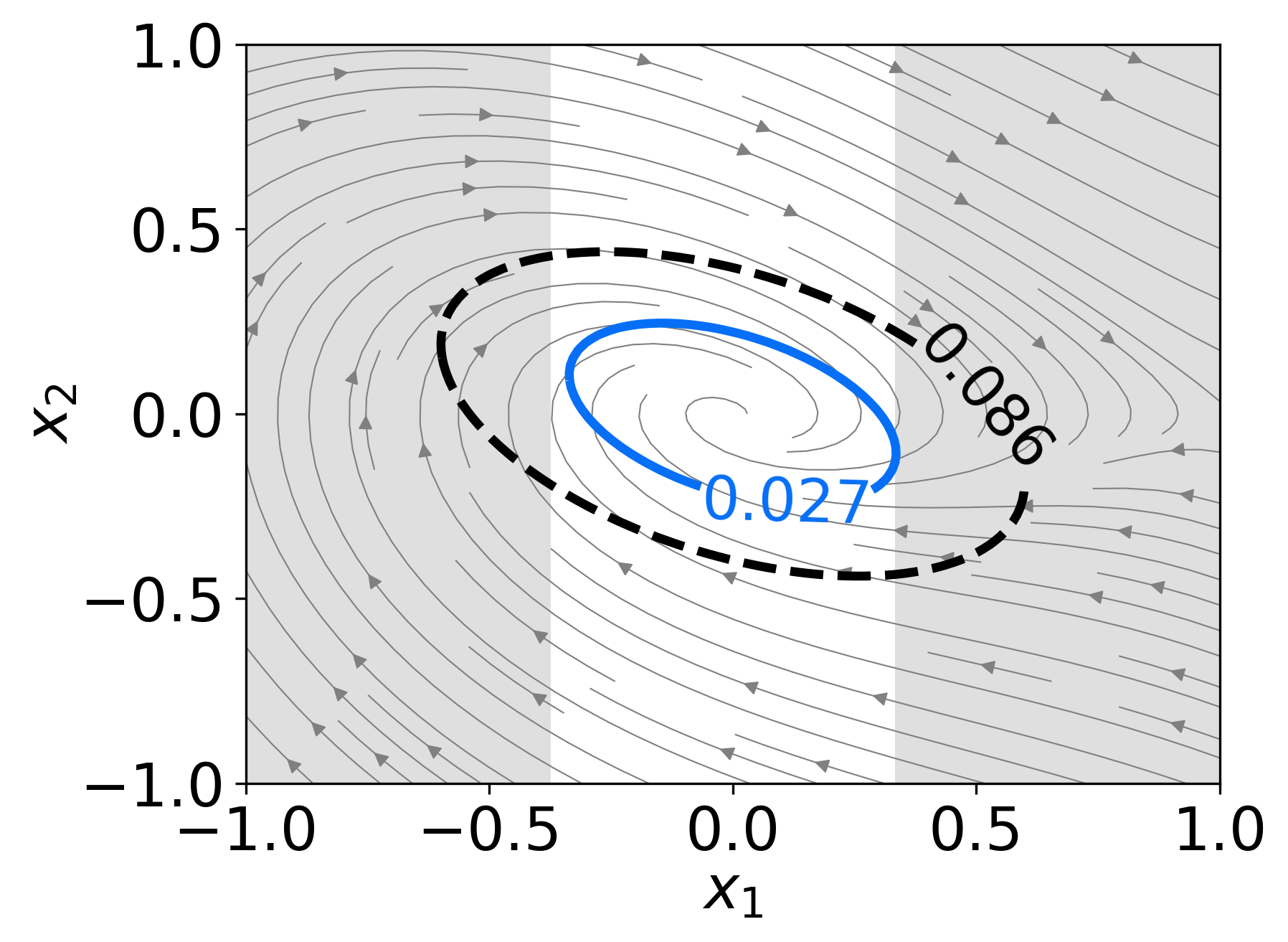}
        (c) Two-machine Power
    \end{minipage}
    \hfill
    \begin{minipage}[t]{0.22\textwidth}
        \centering
        \includegraphics[width=\textwidth]{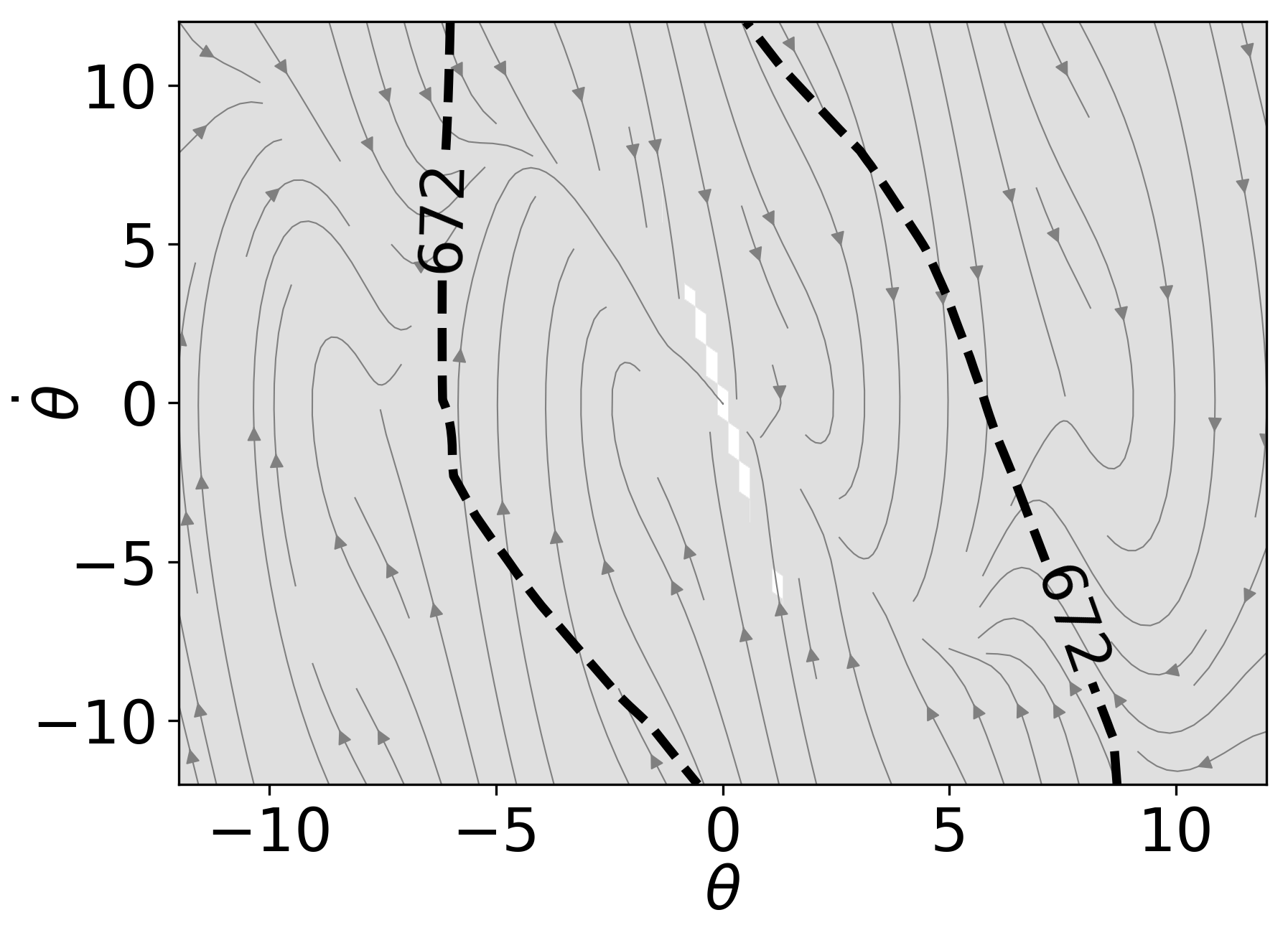}
        (d) Inverted Pendulum
    \end{minipage}
    \vspace{-2mm}
    \caption{The figures visualize the certified region of contraction for our neural contraction metrics and the baseline constant metrics. For each $x \in \mathcal{B}$, we simulate the minimum of $G(x,\delta)$ for a small neighborhood of $\delta$. The grey region indicates $G(x,\delta) < 0$, and the white region indicates the converse. In principle, any forward invariant set inside the white region can be certified as contraction. The top rows visualize the behavior of our learned neural contraction metrics, whereas the bottom row visualizes thebaseline of using a constant contraction metric.}
    \label{fig:all_columns}
\end{figure}

\textbf{Van der Pol Equation.} We consider the Euler discretization of the reverse Van der Pol equation:
\begin{equation}
    \begin{aligned}
        \dot{x_1} = -x_2, \quad \dot{x_2} = x_1 -\mu\left(1 - x_1^2\right)x_2
    \end{aligned}
\end{equation}
with a step size of 0.05. Here we set $\mu = 3$ following \cite{10714396}. On the domain $\mathcal{B} = [-1.2, 1.2] \times [-2.3, 2.3]$, a quadratic Lyapunov function is employed to identify a forward invariant subset of the form $V_{\rho} = \{x \mid V(x) <\rho\}$ with $\rho = 0.119$. A metric that certifies contraction over the entire set $V_{\rho}$ is found. Figure~\ref{fig:all_columns} shows the verified invariant sets where the contraction condition given by the learned NN contraction metric is satisfied (blue ellipse), which occupies the full forward invariant set (dashed black ellipse). In comparison, the largest verified invariant set satisfying the contraction condition given by the constant metric (depicted in the bottom row) is significantly smaller. Figure~\ref{fig:all_columns} further illustrates the regions where the contraction condition given by the learned metric is satisfied (white) and unsatisfied (grey). 
\medskip
\noindent\textbf{Polynomial System.} We consider the Euler discretization of the polynomial system from \cite{10714396}:
\begin{equation}
    \begin{aligned}
        \dot{x}_1 = x_2, \quad \dot{x}_2 = -2x_1 + \frac{1}{3}x_1^3 - x_2
    \end{aligned}
\end{equation}
 with a step size of 0.05. On a domain $\mathcal{B}=[-4, 4]\times [-4, 4]$, we first compute a quadratic Lyapunov function with verified ROA given by sublevel set $V_{\rho}$ with $\rho = 10.8$. With a NN contraction metric, we're able to verify contraction on a sublevel set with $\rho = 10.0$, which occupies over 90 percent of the ROA. In contrast, with a constant metric, the largest sublevel set satisfying the contraction condition is only $\rho=1.5$, as shown in Figure~\ref{fig:all_columns}.

\medskip
\noindent\textbf{Two-machine Power System.} We consider the Euler discretization of the two-machine power system from \cite{10714396}, using a step size of 0.05. The system is described by
\begin{equation}
    \begin{aligned}
        \dot{x_1} = x_2, \quad \dot{x_2} = -0.5x_2 - (\sin(x_1 + \delta) - \sin(\delta))
    \end{aligned}
\end{equation}
with $\delta=\pi/3$. Similar to the examples above, on $\mathcal{B} = [-1, 1]\times [-1,1]$, we identify a quadratic Lyapunov function with verified ROA $V_\rho$ given by $\rho = 0.086$. With a NN contraction metric, we're able to verify contraction on the whole ROA. In comparison with a constant metric, the largest sublevel $V_\rho$ satisfying the contraction is the one with $\rho=0.027$, which occupies only 30 percent of the ROA as shown in Figure~\ref{fig:all_columns}.


\subsection{Non-Smooth Neural Network Controlled System}
\textbf{Inverted Pendulum} We consider the Euler discretization of the controlled inverted pendulum
\begin{align}
\ddot{\theta} &= -\frac{\beta}{ml^2}\dot{\theta} + \frac{g}{l}\sin(\theta) + \frac{u}{ml^2}
\end{align}
with a step size of 0.05. We took the pretrained NN-based state feedback controller and an NN Lyapunov function given in~\cite{yang2024lyapunov}, which achieved state-of-the-art performance in discrete time ROA synthesis. Thanks to the flexibility of our novel formulation and the efficiency of $\alpha,\!\beta$-CROWN, we can also find contraction metrics for NN-controlled state feedback systems, which are often highly non-linear and \emph{non-smooth}. The largest ROA is given by $V_\rho$ with $\rho=672$. Using our pipeline, we successfully learn a contraction metric that certifies the contraction of the system in the sublevel set with $\rho=580$, which occupies around $90\%$ percent of the ROA as shown in Figure~\ref{fig:all_columns}. A constant quadratic metric doesn't provide any reasonable region that satisfies the contraction condition. To the best of our knowledge, this provides the first certified region of contraction for an NN-controlled state feedback nonlinear system.

\section{Conclusion}\label{sec:conclusion}
To summarize, we examine the stability of controlled discrete-time systems using contraction metrics, proposing a novel formulation that certifies contraction under the weaker assumption of continuity in the dynamics. By integrating the advanced neural network verifier $\alpha,\!\beta$-CROWN, we enhance the training and verification process and certify contraction over nearly the entire forward-invariant set. Notably, this approach provides the first formally certified contraction metrics for a discrete-time neural network controlled state-feedback system in the literature.

\section*{Acknowledgments}
B. Hu is generously supported by 
the AFOSR award FA9550-23-1-0732.
H. Zhang is supported in part by the AI2050 program at Schmidt Sciences (AI2050 Early Career Fellowship) and NSF (IIS-2331967).

\bibliography{reference}

\begin{thebibliography}{34}
\providecommand{\natexlab}[1]{#1}
\providecommand{\url}[1]{\texttt{#1}}
\expandafter\ifx\csname urlstyle\endcsname\relax
  \providecommand{\doi}[1]{doi: #1}\else
  \providecommand{\doi}{doi: \begingroup \urlstyle{rm}\Url}\fi

\bibitem[Araujo et~al.(2023)Araujo, Havens, Delattre, Allauzen, and Hu]{araujo2023a}
Alexandre Araujo, Aaron~J Havens, Blaise Delattre, Alexandre Allauzen, and Bin Hu.
\newblock A unified algebraic perspective on {L}ipschitz neural networks.
\newblock In \emph{International Conference on Learning Representations}, 2023.

\bibitem[Bullo(2022)]{bullo2022contraction}
Francesco Bullo.
\newblock \emph{Contraction theory for dynamical systems}.
\newblock Francesco Bullo, 2022.

\bibitem[Chang et~al.(2019)Chang, Roohi, and Gao]{chang2019neural}
Ya-Chien Chang, Nima Roohi, and Sicun Gao.
\newblock Neural lyapunov control.
\newblock \emph{Advances in neural information processing systems}, 32, 2019.

\bibitem[Dai et~al.(2021)Dai, Landry, Yang, Pavone, and Tedrake]{dai2021lyapunov}
Hongkai Dai, Benoit Landry, Lujie Yang, Marco Pavone, and Russ Tedrake.
\newblock Lyapunov-stable neural-network control.
\newblock \emph{arXiv preprint arXiv:2109.14152}, 2021.

\bibitem[Dawson et~al.(2023)Dawson, Gao, and Fan]{dawson2023safe}
Charles Dawson, Sicun Gao, and Chuchu Fan.
\newblock Safe control with learned certificates: A survey of neural lyapunov, barrier, and contraction methods for robotics and control.
\newblock \emph{IEEE Transactions on Robotics}, 39\penalty0 (3):\penalty0 1749--1767, 2023.

\bibitem[Everett et~al.(2023)Everett, Bunel, and Omidshafiei]{everett2023drip}
Michael Everett, Rudy Bunel, and Shayegan Omidshafiei.
\newblock Drip: Domain refinement iteration with polytopes for backward reachability analysis of neural feedback loops.
\newblock \emph{IEEE Control Systems Letters}, 7:\penalty0 1622--1627, 2023.

\bibitem[Fazlyab et~al.(2019)Fazlyab, Robey, Hassani, Morari, and Pappas]{fazlyab2019efficient}
Mahyar Fazlyab, Alexander Robey, Hamed Hassani, Manfred Morari, and George Pappas.
\newblock Efficient and accurate estimation of {L}ipschitz constants for deep neural networks.
\newblock \emph{Advances in neural information processing systems}, 32, 2019.

\bibitem[Fitzsimmons and Liu(2024)]{10714396}
Maxwell Fitzsimmons and Jun Liu.
\newblock Computation and formal verification of neural network contraction metrics.
\newblock \emph{IEEE Control Systems Letters}, pages 1--1, 2024.
\newblock \doi{10.1109/LCSYS.2024.3478272}.

\bibitem[Giesl and Wendland(2019)]{giesl2019construction}
Peter Giesl and Holger Wendland.
\newblock Construction of a contraction metric by meshless collocation.
\newblock \emph{Discrete and Continuous Dynamical Systems-B}, 24\penalty0 (8):\penalty0 3843--3863, 2019.

\bibitem[Giesl et~al.(2022)Giesl, Hafstein, and Kawan]{giesl2022review}
Peter Giesl, Sigurdur Hafstein, and Christoph Kawan.
\newblock Review on contraction analysis and computation of contraction metrics.
\newblock \emph{arXiv preprint arXiv:2203.01367}, 2022.

\bibitem[Giesl et~al.(2024)Giesl, Hafstein, and Mehrabinezhad]{giesl2024contraction}
Peter Giesl, Sigurdur Hafstein, and Iman Mehrabinezhad.
\newblock Contraction metric computation using numerical integration and quadrature.
\newblock \emph{Discrete and Continuous Dynamical Systems-B}, 29\penalty0 (6):\penalty0 2610--2632, 2024.

\bibitem[Huang et~al.(2021)Huang, Zhang, Shi, Kolter, and Anandkumar]{huang2021training}
Yujia Huang, Huan Zhang, Yuanyuan Shi, J~Zico Kolter, and Anima Anandkumar.
\newblock Training certifiably robust neural networks with efficient local {L}ipschitz bounds.
\newblock \emph{Advances in Neural Information Processing Systems}, 34:\penalty0 22745--22757, 2021.

\bibitem[Lee(2018)]{lee2018introduction}
John~M Lee.
\newblock \emph{Introduction to Riemannian manifolds}, volume~2.
\newblock Springer, 2018.

\bibitem[Lohmiller and Slotine(1998)]{lohmiller1998contraction}
Winfried Lohmiller and Jean-Jacques~E Slotine.
\newblock On contraction analysis for non-linear systems.
\newblock \emph{Automatica}, 34\penalty0 (6):\penalty0 683--696, 1998.

\bibitem[Lyapunov(1992)]{lyapunov1992general}
Aleksandr~Mikhailovich Lyapunov.
\newblock The general problem of the stability of motion.
\newblock \emph{International journal of control}, 55\penalty0 (3):\penalty0 531--534, 1992.

\bibitem[Manchester and Slotine(2017)]{manchester2017control}
Ian~R Manchester and Jean-Jacques~E Slotine.
\newblock Control contraction metrics: Convex and intrinsic criteria for nonlinear feedback design.
\newblock \emph{IEEE Transactions on Automatic Control}, 62\penalty0 (6):\penalty0 3046--3053, 2017.

\bibitem[Pokkakkillath and Giesl(2024)]{pokkakkillath2024construction}
Sareena Pokkakkillath and Peter Giesl.
\newblock Construction of contraction metrics for discrete-time dynamical systems using meshfree collocation.
\newblock \emph{Discrete and Continuous Dynamical Systems-B}, 29\penalty0 (4):\penalty0 2043--2071, 2024.

\bibitem[Simpson-Porco and Bullo(2014)]{simpson2014contraction}
John~W Simpson-Porco and Francesco Bullo.
\newblock Contraction theory on riemannian manifolds.
\newblock \emph{Systems \& Control Letters}, 65:\penalty0 74--80, 2014.

\bibitem[Sun et~al.(2021)Sun, Jha, and Fan]{sun2021learning}
Dawei Sun, Susmit Jha, and Chuchu Fan.
\newblock Learning certified control using contraction metric.
\newblock In \emph{Conference on Robot Learning}, pages 1519--1539. PMLR, 2021.

\bibitem[Tran et~al.(2018)Tran, R{\"u}ffer, and Kellett]{tran2018convergence}
Duc~N Tran, Bj{\"o}rn~S R{\"u}ffer, and Christopher~M Kellett.
\newblock Convergence properties for discrete-time nonlinear systems.
\newblock \emph{IEEE Transactions on Automatic Control}, 64\penalty0 (8):\penalty0 3415--3422, 2018.

\bibitem[Tsukamoto and Chung(2020)]{tsukamoto2020neural}
Hiroyasu Tsukamoto and Soon-Jo Chung.
\newblock Neural contraction metrics for robust estimation and control: A convex optimization approach.
\newblock \emph{IEEE Control Systems Letters}, 5\penalty0 (1):\penalty0 211--216, 2020.

\bibitem[Tsukamoto et~al.(2021)Tsukamoto, Chung, and Slotine]{tsukamoto2021contraction}
Hiroyasu Tsukamoto, Soon-Jo Chung, and Jean-Jaques~E Slotine.
\newblock Contraction theory for nonlinear stability analysis and learning-based control: A tutorial overview.
\newblock \emph{Annual Reviews in Control}, 52:\penalty0 135--169, 2021.

\bibitem[Wang and Fazlyab(2024)]{wang2024actor}
Jiarui Wang and Mahyar Fazlyab.
\newblock Actor-critic physics-informed neural lyapunov control.
\newblock \emph{arXiv preprint arXiv:2403.08448}, 2024.

\bibitem[Wang et~al.(2021)Wang, Zhang, Xu, Lin, Jana, Hsieh, and Kolter]{wang2021beta}
Shiqi Wang, Huan Zhang, Kaidi Xu, Xue Lin, Suman Jana, Cho-Jui Hsieh, and J~Zico Kolter.
\newblock {Beta-CROWN}: Efficient bound propagation with per-neuron split constraints for complete and incomplete neural network verification.
\newblock \emph{Advances in Neural Information Processing Systems}, 34, 2021.

\bibitem[Wang et~al.(2024)Wang, Hu, Havens, Araujo, Zheng, Chen, and Jha]{wang2024scalability}
Zi~Wang, Bin Hu, Aaron~J Havens, Alexandre Araujo, Yang Zheng, Yudong Chen, and Somesh Jha.
\newblock On the scalability and memory efficiency of semidefinite programs for {L}ipschitz constant estimation of neural networks.
\newblock In \emph{International Conference on Learning Representations}, 2024.

\bibitem[Wei et~al.(2021)Wei, Mccloy, and Bao]{wei2021control}
Lai Wei, Ryan Mccloy, and Jie Bao.
\newblock Control contraction metric synthesis for discrete-time nonlinear systems.
\newblock \emph{IFAC-PapersOnLine}, 54\penalty0 (3):\penalty0 661--666, 2021.

\bibitem[Wei et~al.(2022)Wei, McCloy, and Bao]{wei2022discrete}
Lai Wei, Ryan McCloy, and Jie Bao.
\newblock Discrete-time contraction-based control of nonlinear systems with parametric uncertainties using neural networks.
\newblock \emph{Computers \& Chemical Engineering}, 166:\penalty0 107962, 2022.

\bibitem[Wu et~al.(2023)Wu, Clark, Kantaros, and Vorobeychik]{wu2023neural}
Junlin Wu, Andrew Clark, Yiannis Kantaros, and Yevgeniy Vorobeychik.
\newblock Neural lyapunov control for discrete-time systems.
\newblock \emph{Advances in neural information processing systems}, 36:\penalty0 2939--2955, 2023.

\bibitem[Xu et~al.(2020)Xu, Shi, Zhang, Wang, Chang, Huang, Kailkhura, Lin, and Hsieh]{xu2020automatic}
Kaidi Xu, Zhouxing Shi, Huan Zhang, Yihan Wang, Kai-Wei Chang, Minlie Huang, Bhavya Kailkhura, Xue Lin, and Cho-Jui Hsieh.
\newblock Automatic perturbation analysis for scalable certified robustness and beyond.
\newblock \emph{Advances in Neural Information Processing Systems}, 33, 2020.

\bibitem[Xu et~al.(2021)Xu, Zhang, Wang, Wang, Jana, Lin, and Hsieh]{xu2021fast}
Kaidi Xu, Huan Zhang, Shiqi Wang, Yihan Wang, Suman Jana, Xue Lin, and Cho-Jui Hsieh.
\newblock {Fast and Complete}: Enabling complete neural network verification with rapid and massively parallel incomplete verifiers.
\newblock In \emph{International Conference on Learning Representations}, 2021.

\bibitem[Yang et~al.(2024)Yang, Dai, Shi, Hsieh, Tedrake, and Zhang]{yang2024lyapunov}
Lujie Yang, Hongkai Dai, Zhouxing Shi, Cho-Jui Hsieh, Russ Tedrake, and Huan Zhang.
\newblock Lyapunov-stable neural control for state and output feedback: A novel formulation.
\newblock In \emph{Forty-first International Conference on Machine Learning}, 2024.

\bibitem[Yang et~al.(2013)Yang, Liu, and Huang]{yang2013neural}
Xiong Yang, Derong Liu, and Yuzhu Huang.
\newblock Neural-network-based online optimal control for uncertain non-linear continuous-time systems with control constraints.
\newblock \emph{IET Control Theory \& Applications}, 7\penalty0 (17):\penalty0 2037--2047, 2013.

\bibitem[Zhang et~al.(2018)Zhang, Weng, Chen, Hsieh, and Daniel]{zhang2018efficient}
Huan Zhang, Tsui-Wei Weng, Pin-Yu Chen, Cho-Jui Hsieh, and Luca Daniel.
\newblock Efficient neural network robustness certification with general activation functions.
\newblock \emph{Advances in Neural Information Processing Systems}, 31:\penalty0 4939--4948, 2018.

\bibitem[Zhang et~al.(2022)Zhang, Wang, Xu, Li, Li, Jana, Hsieh, and Kolter]{zhang2022general}
Huan Zhang, Shiqi Wang, Kaidi Xu, Linyi Li, Bo~Li, Suman Jana, Cho-Jui Hsieh, and J~Zico Kolter.
\newblock General cutting planes for bound-propagation-based neural network verification.
\newblock \emph{Advances in Neural Information Processing Systems}, 2022.

\end{thebibliography}

\end{document}